\documentclass[11pt]{article}
\usepackage{amsfonts, amssymb, latexsym}

\usepackage{color}

\usepackage{color}
\usepackage{mathrsfs}

\def\R{\mathbb{R}}

\def\N{\mathbb{N}}

\def\Z2{{{\mathbb{Z}}^2}}

\def\Z{\mathbb{Z}}
\def\bZ{{\mathbf Z}}

\def\bC{{\mathbf C}}

\def\S{\mathbb{S}}
\def\0{{\bf 0}}
\def\x{{\bf x}}

\def\rE{\mathrm{E}}     %expectation,

\def\v{\vartheta}

\def\R{\mathbb{R}}

\def\N{\mathbb{N}}

\def\Z2{{{\mathbb{Z}}^2}}

\def\Z{\mathbb{Z}}
\def\bZ{{\mathbf Z}}

\def\bC{{\mathbf C}}

\def\S{\mathbb{S}}
\def\0{{\bf 0}}
\def\x{{\bf x}}

\def\rE{\mathrm{E}}     %expectation,

\def\v{\vartheta}

\def\R{\mathbb{R}}

\def\N{\mathbb{N}}

\def\Z2{{{\mathbb{Z}}^2}}

\def\Z{\mathbb{Z}}
\def\bZ{{\mathbf Z}}

\def\bC{{\mathbf C}}

\def\u{\mathbf{u}}
\def\v{\mathbf{v}}
\def\S{\mathbb{S}}
\def\T{\mathbb{T}}
\def\0{{\bf 0}}
\def\x{{\bf x}}

\def\w{\mathbf{w}}
\def\rE{\mathrm{E}}     %expectation,

\def\0{{\bf 0}}
\def\rE{\mathrm{E}}     %expectation,
\def\var{\mathop{\rm var}\nolimits}    %variance
\def\cov{\mathop{\rm cov}\nolimits}    %covariance

\usepackage{amsfonts, amssymb, latexsym}
\usepackage{color}
\usepackage{amsfonts, amssymb, latexsym}

\usepackage{epsfig,graphicx}

\def\0{{\bf 0}}
\def\rE{\mathrm{E}}     %expectation,
\def\var{\mathop{\rm var}\nolimits}    %variance
\def\cov{\mathop{\rm cov}\nolimits}    %covariance

\date{}

\begin{document}

\title{Time varying axially  symmetric vector random fields  \\
         on the sphere }

\author{
   Chunsheng Ma\footnote{  
   Chunsheng  Ma  
     ~~
Department of Mathematics, Statistics, and Physics, Wichita State
University, Wichita, Kansas 67260-0033, USA
~~ e-mail:  chunsheng.ma@wichita.edu 
}
}

\maketitle
\date{}

~ \\

~ \\

\noindent 
{\bf Abstract}    ~~  
This paper presents 
a general  form of the  covariance matrix structure for
 a vector random field  that is axially symmetric   and mean square continuous on  the sphere and
   provides a series representation for a  longitudinally reversible one.
  The  series representation  is somehow  an  imitator of the covariance matrix function, and  both of them have simpler forms  than those proposed 
  in the literature   in terms of  the associated Legendre functions  and  are   useful for modeling and simulation.
  Also,  a general form of the covariance matrix structure 
  is  derived     for
 a  spatio-temporal vector random field   that is axially symmetric   and mean square continuous over  the sphere,
 and   a series representation is given for a   longitudinally reversible one.

~


~

\noindent
{\bf Mathmatics Subject Classification (2010)}   ~~   60G60  $\cdot$ 62M10  $\cdot$  62M30

\newpage

\section{Introduction}

Axially symmetric random fields on the three-dimensional sphere  were introduced by \cite{Jones1963} more than fifty years ago. 
This kind of spatial or spatio-temporal models looks like  reasonable for the practical applications (\cite{Castruccio2013}, \cite{Hitczenko2012},  \cite{Jun2011}, \cite{Jun2007}, \cite{Jun2008},  \cite{Stein2007}), 
 but its complicated  covariance  structure and
series expansion  (\cite{Huang2012}, \cite{Jones1963},  \cite{Stein2007}) have made it very hard to fit real data.
This calls for a more detailed investigation,  and motivates us here to search for    simple forms of the correlation structure and
series expansion  for a scalar, vector, or time varying random field  that is axially symmetric and  mean square continuous
 on the  sphere. 
 The established theories differ from, but are not more complicated than,  those of the  scalar, vector, or time varying random field  that is
  isotropic  and  mean square continuous
 on spheres (\cite{Askey1976}, \cite{Bingham1973}, \cite{Cheng2016},  \cite{Cohen2012},   \cite{Dovidio2014},  \cite{DuMaLi2013},  \cite{Gaspari1999},
 \cite{Gaspari2006}, \cite{Hannan1970},  \cite{Leonenko2012},  \cite{Leonenko2013}, 
 \cite{Ma2012} - %\cite{Ma2016c},  \cite{Malyarenko2013}, 
 \cite{Malyarenko1992}, \cite{Mokljacuk1979},  \cite{Roy1973}, \cite{Roy1976}, \cite{Yadrenko1983}, \cite{Yaglom1961}, \cite{Yaglom1987}).

Denote by $\S^2$ the spherical shell of radius 1 and center $\mathbf{0}$ in $\R^3$, i.e., $\S^2= \{ \| \x \| =1, \x \in \R^3 \}$,
where  $\| \x \|$ is the Euclidean norm of $\x \in \R^3$. 
 Using  spherical coordinates,  a point  $\x \in \S^2$ is determined by the longitude $\theta$ and the latitude $\varphi$, 
 and is designated as $\x = (\varphi, \theta)$,
 where  $0 \le \theta \le 2 \pi$ and
 $0 \le  \varphi \le \pi$.
 
 An $m$-variate real and second-order random field $\{ \bZ (\x), \x = (\varphi, \theta) \in \S^2 \}$ is said to be axially symmetric, if it mean 
 function $\rE \bZ (\x)$ depends only on $\varphi$, and its covariance matrix function $\cov ( \bZ (\x_1), \bZ (\x_2))$ depends on
 $\varphi_1, \varphi_2$, and $\theta_1-\theta_2$.  The covariance matrix function is denoted by $\bC (\varphi_1, \varphi_2, \theta)$, or
   $$ \bC (\varphi_1, \varphi_2,  \theta_1-\theta_2) = \rE \{  (\bZ (\x_1) - \rE  \bZ (\x_1))  (\bZ (\x_2) - \rE  \bZ (\x_2))' \} ,   ~~ \x_1 = (\varphi_1,  \theta_1),  \x_2 = (\varphi_2, \theta_2)  \in \S^2. $$
   Moreover, the random field or  its covariance matrix function is said to be longitudinally reversible \cite{Stein2007},  if 
   $$  \bC( \varphi_1, \varphi_2, -\theta) = \bC (\varphi_1, \varphi_2, \theta),     ~~~~~~
        \varphi_1, \varphi_2 \in [0, \pi],   ~ \theta \in [- 2 \pi,  2 \pi]. $$  
        
In a spatio-temporal setting,        an $m$-variate real and second-order random field $\{ \bZ (\x; t), \x  \in \S^2, t \in \T \}$ is said to be axially symmetric on $\S^2$ and stationary over the temporal domain $\T$, if it mean 
 function $\rE \bZ (\x; t)$ depends only on $\varphi$, and its covariance matrix function $\cov ( \bZ (\x_1; t_1), \bZ (\x_2; t_2))$ depends on
 $\varphi_1, \varphi_2$,  $\theta_1-\theta_2$, and $t_1-t_2$, where $\T = \R$ or $\T$.  The latter is denoted by $\bC (\varphi_1, \varphi_2, \theta; t)$, or
   $$ \bC (\varphi_1, \varphi_2,  \theta_1-\theta_2; t_1 -t_2) = \rE \{  (\bZ (\x_1; t_1) - \rE  \bZ (\x_1; t_1))  (\bZ (\x_2; t_2) - \rE  \bZ (\x_2; t_2))' \} ,    $$
    \hfill  ~~ $\x_k = (\varphi_k,  \theta_k)   \in \S^2, ~  t_k \in \T,  ~ k =1, 2. $
    
    \noindent
  Two fundamental properties of this $m \times m$ matrix function are: (i) $\bC ( \varphi_1, \varphi_2, \theta; t) = ( \bC (\varphi_2, \varphi_1, -\theta; -t) )'$, and (ii) inequality
    \begin{equation}
     \label{positive.definite}
      \sum_{i=1}^n \sum_{j=1}^n  \mathbf{a}'_i \bC ( \varphi_i, \varphi_j, \theta_i-\theta_j; t_i -t_j) \mathbf{a}_j \ge 0  
     \end{equation}
      holds for every   $n \in \N$,  any  $\varphi_i \in [0, \pi], \theta_i \in [0, 2 \pi],$ $t_i \in \T$,  and $\mathbf{a}_i \in \R^m$ ($ i =1, 2, \ldots, n$), where $\N$ stands for the set of positive integers. 
      On the other hand, given an $m \times m$ matrix function with these properties, there exists an $m$-variate Gaussian or elliptically contoured random field   $\{ \bZ (\x;  t), \x \in \S^2, t \in \T \}$
      with  $\bC ( \varphi_1, \varphi_2, \theta; t)$ as its covariance matrix function \cite{Ma2011}.

  Section 2  derives a general form of the covariance matrix structure of an $m$-variate  axially symmetric and mean square continuous random field on $\S^2$ and a series representation of a longitudinally reversible one, starting from a close look at
   the scalar case of \cite{Jones1963}. 
  Moving successively, a general covariance matrix form is given for 
    an $m$-vatiate spatio-temporal random field that is axially symmetric and mean square continuous on $\S^2$ and stationary on $\T$, and a series representation is offered to a longitudinally reversible one over the $\S^2$.
    Theorems 1-3 and 5-7 are  proved in Section 4.

\section{Covariance  matrix structures and series representations}

  Suppose that  $\{ \bZ (\x), \x  \in \S^2 \}$ is an $m$-variate axially symmetric and mean square continuous random field. 
  This section derives a general form of its covariance matrix function, which does not look as complicated as that in \cite{Jones1963} does,  and presents a series representation for a  longitudinally reversible one, which  is useful for modeling and simulation.  An isotropic version is treated in \cite{Ma2012}, \cite{Ma2016a}.

  ~
  
  To begin with, let us take a look at (12) of \cite{Jones1963} for the covariance function of a scalar axially symmetric random field
  $\{ Z(\varphi, \theta), \varphi \in [0, \pi], \theta \in [0, 2 \pi] \} $, which reads  (in our notions)
      \begin{equation}
 \label{Jones1963.form}
  C( \varphi_1, \varphi_2, \theta) = 
  \sum_{n=0}^\infty ( b_n (\varphi_1, \varphi_2) \cos (n \theta)
      +a_n   (\varphi_1, \varphi_2) \sin (n \theta) ),  ~~ \varphi_1, \varphi_2 \in [0, \pi], \theta \in [-2 \pi, 2 \pi],
  \end{equation}
  where 
    \begin{eqnarray*}
     b_n (\varphi_1, \varphi_2) & = & \sum_{i=n}^\infty \sum_{j=n}^\infty b_{n,ij} P_i^n (\cos \varphi_1)  P_j^n (\cos \varphi_2), \\
     a_n (\varphi_1, \varphi_2)  & = & \sum_{i=n}^\infty \sum_{j=n}^\infty a_{n, ij} P_i^n (\cos \varphi_1)  P_j^n (\cos \varphi_2), ~~~~ n \in \N_0, 
    \end{eqnarray*} 
    $b_{n, ij}$ and  $a_{n, ij}$ are certain constants,   $P^n_k (x)$ are the associated Legendre function \cite{Andrews1999}, and $\N_0$ denotes the set of nonnegative integers. 
    Simply speaking, (\ref{Jones1963.form}) is a Fourier series expansion of $C(\varphi_1, \varphi_2, \theta)$ in terms of $\theta \in [- 2 \pi, 2 \pi],$  but using or not  the associated Legendre function is not so important there, as is indicated   by Theorem 1.

    ~

    \noindent
    {\bf Theorem 1}  ~~   {\em    For   an  axially symmetric and  mean square continuous random field  $\{  Z (\x), \x \in \S^2 \}$,
     its  covariance function $C(\varphi_1, \varphi_2, \theta)$ is 
      of the form   }
       \begin{equation}
       \label{cov.matrix.fun.1}
          C(\varphi_1, \varphi_2, \theta)   =
            \sum\limits_{n=0}^\infty  \{ b_n (\varphi_1, \varphi_2) \cos (n  \theta)
                           +a_n (\varphi_1, \varphi_2) \sin (n \theta) \},  ~~~   \varphi_1, \varphi_2 \in [0, \pi], ~ \theta \in [- 2 \pi,  2 \pi], 
       \end{equation}
   {\em    where    $ b_n (\varphi_1, \varphi_2)$ is a  covariance  function on $[0, \pi]$  for each fixed $n \in \N_0$, and $\sum\limits_{n=0}^\infty  b_n ( \varphi_1, \varphi_2)$  converges.  }

   ~

     Although  $b_n (\varphi_1, \varphi_2)$ has to be  a covariance function for each $n \in \N_0$,  it is not clear
    how to  interpret $a_n (\varphi_1, \varphi_2)$. 
    No matter what  $a_n (\varphi_1, \varphi_2)$ looks like, it must make (\ref{cov.matrix.fun.1}) a positive definite function.  An example of (\ref{cov.matrix.fun.1}) is illustrated in Example 1,  where $a_n (\varphi_1, \varphi_2)$
    is either non-negatively or negatively proportional to $b_n (\varphi_1, \varphi_2)$.
    
    ~
    
    \noindent
    {\em Example 1}    ~~  For a constant $\lambda$ with $|\lambda| \le 1$,
          $$ C( \varphi_1, \varphi_2, \theta) 
               = \sum\limits_{n=0}^\infty b_n (\varphi_1, \varphi_2) 
                  ( \cos (n \theta)+\lambda \sin (n \theta)),  ~~~ \varphi_1, \varphi_2 \in [0, \pi], ~ \theta \in [- 2 \pi, 2 \pi], $$
          is  the covariance function of an axially symmetric Gaussian or elliptically contoured random field on $\S^2$,         
    where    $ b_n (\varphi_1, \varphi_2)$ is a  covariance  function on $[0, \pi]$  for each fixed $n \in \N_0$, and $\sum\limits_{n=0}^\infty  b_n ( \varphi_1, \varphi_2)$  converges. This follows from Theorem 8 of \cite{Ma2011}, noticing that
    $\cos (n \theta)+\lambda \sin (n \theta)$ is a positive definite function of $\theta \in [- 2 \pi, 2 \pi]$ for each $n \in \N_0$.
    Of course, many   $ b_n (\varphi_1, \varphi_2)$ are available for selection.

  ~
  
  \noindent
  {\bf Corollary 1.1}  ~~ {\em  A  longitudinally reversible  covariance function $C(\varphi_1, \varphi_2, \theta)$  takes the form}
       \begin{equation}
       \label{cov.matrix.fun.1.0}
          C(\varphi_1, \varphi_2, \theta)   =
            \sum\limits_{n=0}^\infty   b_n (\varphi_1, \varphi_2) \cos (n  \theta),  ~~~   \varphi_1, \varphi_2 \in [0, \pi], ~ \theta \in [- 2 \pi,  2 \pi]. 
       \end{equation}

        {\em Conversely, given a function $C(\varphi_1, \varphi_2, \theta)$ of the form (\ref{cov.matrix.fun.1.0}) with 
      summable  $\sum\limits_{n=0}^\infty  b_n ( \varphi_1, \varphi_2)$  and  each $ b_n (\varphi_1, \varphi_2)$ being a covariance  function on   $[0, \pi]$, there exists a  longitudinally reversible
         Gaussian or elliptical contoured
       random field on $\S^2$ with it as the covariance  function.}

       ~

      In fact,  (\ref{cov.matrix.fun.1.0}) follows from (\ref{cov.matrix.fun.1}) and $C(\varphi_1, \varphi_2, \theta)   = C( \varphi_1, \varphi_2, -\theta)$; see also Proposition 3 of \cite{Huang2012}.   
      One way to establish  the second part of  Corollary 1.1 is to verify the  positive definiteness of (\ref{cov.matrix.fun.1.0}) and to use Theorem 8 of \cite{Ma2011}. Instead, we give a series representation for a random field with  (\ref{cov.matrix.fun.1.0}) as its covariance
      function, which is useful for modeling and simulation.
      
      ~

      \noindent
      {\bf Theorem 2}.  {\em   Assume that  stochastic processes  $\{ V_{n1} (\varphi), \varphi \in [0, \pi] \}$ and $\{ V_{n2} (\varphi), \varphi \in [0, \pi] \}$ have mean 0 and covariance function $b_n (\varphi_1, \varphi_2)$  for each $n \in \N_0$, and  that for all $n \in 
      \N_0$,  $\{ V_{n1} (\varphi), \varphi \in [0, \pi] \}$ and $\{ V_{n2} (\varphi), \varphi \in [0, \pi] \}$ are independent.
          If    $\sum\limits_{n=0}^\infty  b_n ( \varphi_1, \varphi_2)$  converges, then}
      \begin{equation}
      \label{series.exp1}
      Z(\varphi, \theta) =  \sum_{n=0}^\infty  ( V_{n1}  ( \varphi) \cos (n \theta) + V_{n2}  ( \varphi) \sin ( n \theta) ),    ~~~~~~  \varphi \in [0, \pi], ~ \theta \in [0, 2 \pi],
       \end{equation}
    {\em is a longitudinally reversible random field on $\S^2$ with mean 0 and covariance function (\ref{cov.matrix.fun.1.0}).}

   ~
   
   Somehow  (\ref{series.exp1}) is an imitator of (\ref{cov.matrix.fun.1}) but not (\ref{cov.matrix.fun.1.0}). It would be of interest to generate an axially symmetric random field on $\S^2$.  The vector case is slightly more complicated than the scalar case, since a covariance matrix function is not necessarily symmetric.

~

    \noindent
    {\bf Theorem 3}  ~~   {\em     If  an  $m$-variate  mean square continuous random field  $\{ \bZ (\x), \x \in \S^2 \}$ is   
     axially symmetric,  then   $\frac{\bC(\varphi_1, \varphi_2, \theta)+ (\bC(\varphi_1, \varphi_2, \theta))'}{2}$ is 
      of the form   }
       \begin{equation}
       \label{cov.matrix.fun.m}
      \begin{array}{c}
       \frac{\bC(\varphi_1, \varphi_2, \theta)+ (\bC(\varphi_1, \varphi_2, \theta))'}{2} 
          =
             \sum\limits_{n=0}^\infty  \left\{   \mathbf{B}_n (\varphi_1, \varphi_2) \cos (n  \theta)
                + \mathbf{A}_n (\varphi_1, \varphi_2) \sin (n  \theta) \right\},   \\
               ~~~~~~~~~~~~~~~~~~~~~~~~~~~~~~~  ~~~   \varphi_1, \varphi_2 \in [0, \pi], ~ \theta \in [- 2 \pi,  2 \pi], 
                 \end{array}      
       \end{equation}
      {\em  where   $ \mathbf{B}_n (\varphi_1, \varphi_2)$ ($ n \in \N_0$) are $m \times m $ symmetric matrices,   $\sum\limits_{n=0}^\infty  \mathbf{B}_n ( \varphi_1, \varphi_2)$  converges,  
     and, for each fixed $n \in \N_0$,  $ \mathbf{B}_n (\varphi_1, \varphi_2)$ is a  covariance matrix function on $[0, \pi]$.  }
   
   ~

    \noindent
    {\em Example 2}    ~~  For a constant $\lambda$ with $|\lambda| \le 1$,
          $$  \bC( \varphi_1, \varphi_2, \theta) 
               = \sum\limits_{n=0}^\infty  \mathbf{B}_n (\varphi_1, \varphi_2) 
                  ( \cos (n \theta)+\lambda \sin (n \theta)),  ~~~ \varphi_1, \varphi_2 \in [0, \pi], ~ \theta \in [- 2 \pi, 2 \pi], $$
          is  the covariance matrix function of an $m$-variate axially symmetric Gaussian or elliptically contoured random field on $\S^2$,         
    if    $ \mathbf{B}_n (\varphi_1, \varphi_2)$ is an $m \times m$  covariance matrix function on $[0, \pi]$  for each fixed $n \in \N_0$, and $\sum\limits_{n=0}^\infty  \mathbf{B}_n ( \varphi_1, \varphi_2)$  converges. 
   
   ~

     \noindent
     {\bf Corolloary 3.1}  ~~   {\em If $\bC(\varphi_1, \varphi_2, \theta) $ is   longitudinally reversible, then }
       \begin{equation}
       \label{cov.matrix.fun.m.0}
          \frac{\bC(\varphi_1, \varphi_2, \theta)+ (\bC(\varphi_1, \varphi_2, \theta))'}{2}     =
            \sum\limits_{n=0}^\infty   \mathbf{B}_n (\varphi_1, \varphi_2) \cos (n  \theta),  ~~~   \varphi_1, \varphi_2 \in [0, \pi], ~ \theta \in [- 2 \pi,  2 \pi]. 
       \end{equation}

       \noindent
       {\bf Corollary 3.2}  ~~ {\em An $m \times m $ matrix function  
       \begin{equation}
       \label{cov.matrix.fun.m.1}
          \bC(\varphi_1, \varphi_2, \theta)   =
            \sum\limits_{n=0}^\infty   \mathbf{B}_n (\varphi_1, \varphi_2) \cos (n  \theta),  ~~~   \varphi_1, \varphi_2 \in [0, \pi], ~ \theta \in [- 2 \pi,  2 \pi], 
       \end{equation}
       is the covariance matrix function of
       an $m$-variate Gaussian or elliptically contoured random field on $\S^2$ if and only if 
        $ \mathbf{B}_n (\varphi_1, \varphi_2)$ is an $m \times m$  covariance matrix function on $[0, \pi]$  for each fixed $n \in \N_0$, and $\sum\limits_{n=0}^\infty  \mathbf{B}_n ( \varphi_1, \varphi_2)$  converges. }
       
      ~
      
      The matrix $\mathbf{B}_n (\varphi_1, \varphi_2)$ in (\ref{cov.matrix.fun.m.0}) has to be symmetric  in the sense that
      $  \mathbf{B}_n (\varphi_1, \varphi_2) = (\mathbf{B}_n (\varphi_1, \varphi_2))'$. But,  it is not necessarily so in (\ref{cov.matrix.fun.m.1}) or  the next theorem, whose proof is similar to that of Theorem 2 and is thus  omitted.  
      
      ~
      
      \noindent
      {\bf Theorem 4}  ~ 
   {\em   Assume that  $m$-variate stochastic processes  $\{ \mathbf{V}_{n1} (\varphi), \varphi \in [0, \pi] \}$ and $\{ \mathbf{V}_{n2} (\varphi), \varphi \in [0, \pi] \}$ have mean $\0$  and covariance matrix function $\mathbf{B}_n (\varphi_1, \varphi_2)$  for each $n \in \N_0$, and  that for all $n \in 
      \N_0$,  $\{ \mathbf{V}_{n1} (\varphi), \varphi \in [0, \pi] \}$ and $\{ \mathbf{V}_{n2} (\varphi), \varphi \in [0, \pi] \}$ are independent.
          If    $\sum\limits_{n=0}^\infty  \mathbf{B}_n ( \varphi_1, \varphi_2)$  converges, then}
      \begin{equation}
      \label{series.expm}
      \bZ(\varphi, \theta) =  \sum_{n=0}^\infty  ( \mathbf{V}_{n1}  ( \varphi) \cos (n \theta) + \mathbf{V}_{n2}  ( \varphi) \sin ( n \theta) ),    ~~~~~~  \varphi \in [0, \pi], ~ \theta \in [0, 2 \pi],
       \end{equation}
    {\em is an $m$-variate  longitudinally reversible random field on $\S^2$ with mean $\0$ and covariance  matrix function
    (\ref{cov.matrix.fun.m.1}).}

                   ~
                   
        \noindent
{\em Example 3}  ~~ Let 
     $$ \mathbf{B}_{n, ij} (\varphi_1, \varphi_2)
          = \left\{
             \begin{array}{ll}
              \frac{1}{ (b_i (\varphi_1)+ b_j (\varphi_2)) \pi },   ~  &   ~ n = 0, \\
              \frac{2}{ (n^2+ b_i (\varphi_1)+ b_j (\varphi_2)) \pi },  ~   &   ~ n \in \N,  
              ~  \varphi_1, \varphi_2 \in [0, \pi], ~~  i, j = 1, \ldots, m,  
             \end{array}
             \right. $$           
    where $b_i (\varphi)$ ($i =1, \ldots, m$) are positive functions on $[0, \pi]$.  
    It can be verified that each       $\mathbf{B}_n (\varphi_1, \varphi_2)$ is a covariance matrix function on $[0, \pi]$.       
    The direct/cross covariance functions of (\ref{series.expm}) are
        \begin{eqnarray*}
           C_{ij} ( \varphi_1, \varphi_2, \theta)   
         &  =  &   \frac{\cosh \left( ( \pi - |\theta|)  \sqrt{b_i(\varphi_1)+b_j (\varphi_2)} \right)}{
                            \sqrt{b_i(\varphi_1)+b_j (\varphi_2)}   \sinh (
                             \sqrt{b_i(\varphi_1)+b_j (\varphi_2)} \pi )},  \\
           &   &     ~~~~~~~~~~~~~~~~~~~   \varphi_1, \varphi_2 \in [0, \pi], ~ \theta \in [- 2 \pi,  2 \pi],  ~  i, j   = 1, \ldots, m, 
           \end{eqnarray*}    
             based on  the identity (see, e.g., page 577 of \cite{Watson1944})
              $$ \sum_{n=1}^\infty \frac{\cos (n \theta)}{n^2+a^2} 
                   = \frac{ \pi \cosh ((\pi-\theta) a)}{ 2 a \sinh (\pi a)}- \frac{1}{2 a^2},
                    ~~~  a > 0, ~ \theta \in [0, 2 \pi]. $$

                   ~

\section{Time varying  axially symmetric vector random fields }

   This section deals with the covariance matrix structure of an $m$-variate random field $\{ \bZ (\x; t),  \x \in \S^2, t \in \T \}$
   that is axially symmetric and mean square  continuous on $\S^2$ and stationary on $\T$, and presents a series representation
   for a  longitudinally reversible one. An isotropic case is studied in  \cite{Ma2016c}. 
   The connection of the Gaussianity between  the random field and its coefficients in the series representation is explored in
   Theorem 6. Some arguments similar to the purely spatial case are not repeated.
   
   ~

    \noindent
    {\bf Theorem 5}  ~~   {\em     If  an  $m$-variate   random field  $\{ \bZ (\x; t), \x \in \S^2, t \in \T \}$ is   
     axially symmetric and  mean square continuous on $\S^2$ and stationary on $\T$,  then    }
       \begin{equation}
       \label{cov.matrix.fun.m.t}
      \begin{array}{c}
       \frac{\bC(\varphi_1, \varphi_2, \theta; t)+ (\bC(\varphi_1, \varphi_2, \theta; t))' + \bC(\varphi_1, \varphi_2, \theta; -t)+ (\bC(\varphi_1, \varphi_2, \theta; -t))'}{2} \\
          =
             \sum\limits_{n=0}^\infty  \left\{   \mathbf{B}_n (\varphi_1, \varphi_2; t) \cos (n  \theta)
                + \mathbf{A}_n (\varphi_1, \varphi_2; t) \sin (n  \theta) \right\},   \\
               ~~~~~~~~~~~~~~~~~~~~~~~~~~~~~~~  ~~~   \varphi_1, \varphi_2 \in [0, \pi], ~ \theta \in [- 2 \pi,  2 \pi],  t \in \T,
                 \end{array}      
       \end{equation}
      {\em  where   $ \mathbf{B}_n (\varphi_1, \varphi_2; t)$ ($ n \in \N_0$) are $m \times m $ symmetric matrices and  $\sum\limits_{n=0}^\infty  \mathbf{B}_n ( \varphi_1, \varphi_2; t)$  converges,  
     and, for each fixed $n \in \N_0$,  $ \mathbf{B}_n (\varphi_1, \varphi_2; t)$ is a  covariance matrix function on $[0, \pi] \times \T$.  }

~

\noindent
{\bf Corollary 5.1} ~~ {\em If  $\{ \bZ (\x; t), \x \in \S^2, t \in \T \}$ is   
     longitudinally reversible on $\S^2$,  then    }
       \begin{equation}
       \label{cov.matrix.fun.m.t.1}
      \begin{array}{c}
       \frac{\bC(\varphi_1, \varphi_2, \theta; t)+ (\bC(\varphi_1, \varphi_2, \theta; t))' + \bC(\varphi_1, \varphi_2, \theta; -t)+ (\bC(\varphi_1, \varphi_2, \theta; -t))'}{2} 
          =
             \sum\limits_{n=0}^\infty     \mathbf{B}_n (\varphi_1, \varphi_2; t) \cos (n  \theta),   \\
               ~~~~~~~~~~~~~~~~~~~~~~~~~~~~~~~  ~~~   \varphi_1, \varphi_2 \in [0, \pi], ~ \theta \in [- 2 \pi,  2 \pi],  t \in \T.
                 \end{array}      
       \end{equation}
       
       ~
       
      \noindent
{\bf Corollary 5.2}  
~~ {\em An $m \times m$ matrix function}
    \begin{equation}
    \label{cov.matrix.fun.m.t.2}
   \bC (\varphi_1, \varphi_2, \theta; t)
             =
            \sum\limits_{n=0}^\infty   \mathbf{B}_n (\varphi_1, \varphi_2; t) \cos (n  \theta),  ~~~   \varphi_1, \varphi_2 \in [0, \pi], ~ \theta \in [- 2 \pi,  2 \pi],  ~ t \in \T,
    \end{equation}
    {\em is the covariance matrix function of an $m$-variate random field  on $\S^2 \times \T$  that is       longitudinally reversible on $\S^2$ and stationary on $\T$  if and only if 
          $\sum\limits_{n=0}^\infty  \mathbf{B}_n ( \varphi_1, \varphi_2; t)$  converges,  
     and, for each fixed $n \in \N_0$,  $ \mathbf{B}_n (\varphi_1, \varphi_2; t)$ is a  covariance matrix function on $[0, \pi] \times \T$.  }

~

 \noindent
      {\bf Theorem 6}  ~ 
   {\em   Assume that  $m$-variate stochastic processes  $\{ \mathbf{V}_{n1} (\varphi; t), \varphi \in [0, \pi],  t \in \T \}$ and $\{ \mathbf{V}_{n2} (\varphi; t), \varphi \in [0, \pi], t \in \T \}$ have mean $\0$  and covariance matrix function $\mathbf{B}_n (\varphi_1, \varphi_2; t)$  for each $n \in \N_0$, and  that for all $n \in 
      \N_0$,  $\{ \mathbf{V}_{n1} (\varphi; t), \varphi \in [0, \pi],  t \in \T \}$ and $\{ \mathbf{V}_{n2} (\varphi; t), \varphi \in [0, \pi],
      t \in \T \}$ are independent.
          If    $\sum\limits_{n=0}^\infty  \mathbf{B}_n ( \varphi_1, \varphi_2; t)$  converges, then}
      \begin{equation}
      \label{series.expm.t}
      \bZ(\varphi, \theta; t) =  \sum_{n=0}^\infty  ( \mathbf{V}_{n1}  ( \varphi; t) \cos (n \theta) + \mathbf{V}_{n2}  ( \varphi; t) \sin ( n \theta) ),    ~~~~~~  \varphi \in [0, \pi], ~ \theta \in [0, 2 \pi], t \in \T,
       \end{equation}
    {\em is an $m$-variate   random field  longitudinally reversible on $\S^2$ and stationary on $\T$,  with mean $\0$ and covariance  matrix function} 
      $$ \bC (\varphi_1, \varphi_2, \theta; t)
             =
            \sum\limits_{n=0}^\infty   \mathbf{B}_n (\varphi_1, \varphi_2; t) \cos (n  \theta),  ~~~   \varphi_1, \varphi_2 \in [0, \pi], ~ \theta \in [- 2 \pi,  2 \pi],  ~ t \in \T.  $$
            
            {\em Moreover,  $\bZ(\varphi, \theta; t)$ is  Gaussian  if and only if all  $\mathbf{V}_{n1} (\varphi; t)$ and $\mathbf{V}_{n2} (\varphi; t)$ ($n \in \N_0$) are Gaussian.}

~

~

\noindent
{\em Example 4}  ~~  Let $\mathbf{B} (\varphi_1, \varphi_2; t)$ be an $m \times m$ covariance matrix function on $[0, \pi] \times \T$ and all its entries be less than 1 in absolute value. 

  \begin{itemize}
  
   \item[(i)] In Theorem 6  choose
   $$ \mathbf{B}_n  ( \varphi_1, \varphi_2; t) =
      \left\{
       \begin{array}{ll}
        \frac{2}{n}  \left( \mathbf{B} (\varphi_1, \varphi_2; t) \right)^{\circ n},   ~   &   ~  n \in \N, \\
        \0,    ~   &  ~  n =0,   ~ \varphi_1, \varphi_2 \in [0, \pi], ~ t \in \T,
        \end{array}
        \right. $$
        where ${\bf B}^{\circ p}$ denotes the  Hadamard  $p$ power of ${\bf B} =(b_{ij})$, whose entries are $b_{ij}^p$, the $p$ power of $b_{ij}, i, j = 1, \ldots, m$.  The direct/cross covariance  functions of (\ref{series.expm.t}) are
        $$ C_{ij} ( \varphi_1, \varphi_2, \theta; t)   = - \ln (1- 2 b_{ij} (\varphi_1, \varphi_2; t) \cos \theta + b^2_{ij} (\varphi_1, \varphi_2; t) ),  $$
             \hfill  $   ~~~   \varphi_1, \varphi_2 \in [0, \pi], ~ \theta \in [- 2 \pi,  2 \pi],  ~  t \in \T,   i, j   = 1, \ldots, m, $
               
               \noindent
             in view of the identity
              $$ \sum_{n=1}^\infty \frac{ a^n}{n} \cos (n \theta) = -\frac{1}{2} \ln ( 1- 2 a \cos \theta+a^2),
                   ~~~~ | a | < 1. $$
                   
       \item[(ii)]      Choose
   $$ \mathbf{B}_n  ( \varphi_1, \varphi_2) =
      \left\{
       \begin{array}{ll}
          \left( \mathbf{B} (\varphi_1, \varphi_2; t) \right)^{\circ n},   ~   &   ~  n \in \N, \\
        {\bf 1},    ~   &  ~  n =0,   ~ \varphi_1, \varphi_2 \in [0, \pi], ~  t \in \T,
        \end{array}
        \right. $$
        where ${\bf 1}$  is an $m \times m$ matrix with all entries 1.  Then  the  direct/cross covariance  functions of (\ref{series.expm.t}) are
        $$ C_{ij} ( \varphi_1, \varphi_2, \theta; t)   = 
               \frac{1-b^2_{ij} (\varphi_1, \varphi_2; t)}{ 1- 2 b_{ij} (\varphi_1, \varphi_2; t) \cos \theta +
                   b^2_{ij} (\varphi_1, \varphi_2; t)},  $$ 
               \hfill  $  ~~~   \varphi_1, \varphi_2 \in [0, \pi], ~ \theta \in [- 2 \pi,  2 \pi],  t \in \T,  ~ i, j   = 1, \ldots, m, $
               
               \noindent
            which follows from the identity
              $$ \sum_{n=0}^\infty  a^n \cos (n \theta) =  \frac{1-a^2}{1-2 a \cos \theta +a^2},
                   ~~~~ | a | < 1. $$
     \end{itemize}              
                   
                   ~
                   
                    A sufficient condition is given in the following theorem for  an $m \times m$ matrix function to be the
                    covariance matrix function of an $m$-variate Gaussian or elliptically contoured random field on $\S^2 \times \T$
                    that is axially symmetric and mean square  continuous on $\S^2$ and stationary on $\T$.
    
    ~
    
    \noindent
    {\bf Theorem 7}. {\em For each $n \in \N_0$, if the inequality}
         \begin{equation}
         \label{suff.con}
         \sum_{i=1}^l \sum_{j=1}^l ( \mathbf{u}'_i  \mathbf{B}_n (\varphi_i, \varphi_j; t_i-t_j) \mathbf{u}_j + \mathbf{v}'_i   \mathbf{B}_n   (\varphi_i, \varphi_j; t_i-t_j)  \mathbf{v}_j + \mathbf{u}'_i  \mathbf{A}_n (\varphi_i, \varphi_j; t_i-t_j) \mathbf{v}_j 
                -  \mathbf{v}'_i   \mathbf{A}_n (\varphi_i, \varphi_j; t_i-t_j) \mathbf{u}_j  ) \ge 0
         \end{equation}
   {\em  holds for every $l \in \N$, any $\mathbf{u}_i \in \R^m,  \mathbf{v}_i \in \R^m$,  $t_i \in \T$, and $\varphi_i \in [0, \pi]$, then
    there exists an $m$-variate  Gaussian or elliptically contoured random field on $\S^2 \times \T$ with
    the covariance matrix function}
           \begin{equation}
       \label{cov.m.t.suf}
      \begin{array}{c}
       \bC(\varphi_1, \varphi_2, \theta; t)
          =  \sum\limits_{n=0}^\infty  \left\{   \mathbf{B}_n (\varphi_1, \varphi_2; t) \cos (n  \theta)
                + \mathbf{A}_n (\varphi_1, \varphi_2; t) \sin (n  \theta) \right\},   \\
               ~~~~~~~~~~~~~~~~~~~~~~~~~~~~~~~  ~~~   \varphi_1, \varphi_2 \in [0, \pi], ~ \theta \in [- 2 \pi,  2 \pi],  t \in \T.
                 \end{array}      
       \end{equation}
       
       ~
       
 Inequality (\ref{suff.con}) ensures that each level $n$ of (\ref{cov.m.t.suf}), 
 $        \mathbf{B}_n (\varphi_1, \varphi_2; t) \cos (n  \theta)
                + \mathbf{A}_n (\varphi_1, \varphi_2; t) \sin (n  \theta)$,  satisfies inequality (\ref{positive.definite}). It would be of interest to see  whether (\ref{suff.con})  is also a necessary condition.
       
       ~

    %   \noindent
    %   {\em Corollary}   ~~ {\em  For an $m \times m$ symmetric matrix $\Lambda$ with all eigenvalues between -1 and 1, } 
    %      $$  \begin{array}{c}
    %   \bC(\varphi_1, \varphi_2, \theta; t)
    %      =  \sum\limits_{n=0}^\infty     \mathbf{B}_n (\varphi_1, \varphi_2; t)  \circ  (  \mathbf{I}_m \cos (n  \theta)
   %             + \Lambda \sin (n  \theta) ),   \\
   %            ~~~~~~~~~~~~~~~~~~~~~~~~~~~~~~~  ~~~   \varphi_1, \varphi_2 \in [0, \pi], ~ \theta \in [- 2 \pi,  2 \pi],  t \in \T,
   %              \end{array} $$
   % {\em  is the covariance matrix function of  an $m$-variate Gaussian or elliptically contoured random field on $\S^2 \times \T$,
  % if   $\sum\limits_{n=0}^\infty  \mathbf{B}_n ( \varphi_1, \varphi_2; t)$  converges,  
  %   and, for each fixed $n \in \N_0$,  $ \mathbf{B}_n (\varphi_1, \varphi_2; t)$ is a  covariance matrix function on $[0, \pi] \times \T$, where $\mathbf{I}_m$ is an identity matrix. }

     ~
     
   %  To verify inequality (\ref{suff.con}), notice that  $\mathbf{A}_n (\varphi_1, \varphi_2; t) = \mathbf{B}_n (\varphi_1, \varphi_2; t) 
   %   \Lambda$ in this case.  

\section{Proofs}

\subsection{Proof of Theorem 1}

Recall that the mean square continuity  of $\{ Z (\x), \x \in \S^2 \}$   means 
    $$ \rE | Z (\varphi_1, \theta_1)- Z(\varphi_2, \theta_2) |^2  \to 0,   ~~~~ \mbox{as}  ~ \varphi_1 \to \varphi_2, \theta_1 \to \theta_2. $$
   This implies that $ C( \varphi_1, \varphi_2, \theta)$ is  continuous with respect to $\theta \in [- 2 \pi, 2 \pi]$
    for fixed $\varphi_1, \varphi_2 \in [0, \pi]$. Indeed, it follows from the  Cauchy-Schwartz inequality that, for $\theta_1, \theta_2
    \ge 0$,
        \begin{eqnarray*}
          &   & | C( \varphi_1, \varphi_2, \theta_1) - C(\varphi_1, \varphi_2, \theta_2)  |   \\
          & =  & | \rE  Z(\varphi_1, \theta_1) (Z(\varphi_2, 0)- \rE  Z(\varphi_2, 0) ) - 
                        \rE  Z(\varphi_1, \theta_2) (Z(\varphi_2, 0)- \rE  Z(\varphi_2, 0) ) | \\
          & =  &  | \rE  \{ (    Z(\varphi_1, \theta_1) -   Z(\varphi_1, \theta_2) )  (Z(\varphi_2, 0)- \rE  Z(\varphi_2, 0) ) \} | \\
          & \le  &    \var (Z (\varphi_2)) \rE |    Z(\varphi_1, \theta_1) -   Z(\varphi_1, \theta_2) |^2 \\
          &  \to  & 0,    ~~~~~~ \theta_1 \to \theta_2,      
        \end{eqnarray*}    
        and,  for $\theta_1, \theta_2 \le 0$,
        \begin{eqnarray*}
          &   & | C( \varphi_1, \varphi_2, \theta_1) - C(\varphi_1, \varphi_2, \theta_2)  |   \\
          & =  & | \rE  (Z(\varphi_1, 0) - \rE  Z(\varphi_1, 0) ) Z(\varphi_2,  -\theta_1) 
                        - \rE ( Z(\varphi_1, 0)- \rE Z(\varphi_1, 0))  Z(\varphi_2, -\theta_2 ) | \\
          & =  &  | \rE \{  (    Z(\varphi_1, 0) -  \rE Z(\varphi_1, 0) )  (Z(\varphi_2, -\theta_1)-  Z(\varphi_2, -\theta_2) ) \} | \\
          & \le  &    \var (Z (\varphi_1)) \rE |    Z(\varphi_2, -\theta_1) -   Z(\varphi_2, -\theta_2) |^2 \\
          &  \to  & 0,    ~~~~~~ \theta_1 \to \theta_2.     
        \end{eqnarray*}    
    Thus,  as a continuous function of $\theta \in [- 2 \pi,  2 \pi]$,  $C(\varphi_1, \varphi_2, \theta)$ takes the Fourier series expansion
    (\ref{cov.matrix.fun.1}), since $\{  \cos (n \theta), \sin (n \theta), n  \in \N_0 \}$ consists of an orthogonal basis of $L^2 [ - 2 \pi,
    2 \pi]$. 
  The convergence of $\sum\limits_{n=0}^\infty b_n (\varphi_1, \varphi_2)$ is due to the existence of $C( \varphi_1, \varphi_2, 0)$.

      What remains is to confirm that each coefficient in (\ref{cov.matrix.fun.1}),  $b_n (\varphi_1, \varphi_2)$,
       is a covariance function.  For this purpose,  consider a  stochastic process
          $$ W_n (\varphi) =  \int_0^{ 2 \pi}  Z (\varphi, \theta)  \cos (n \theta)  d \theta,    ~~~~~~ \varphi \in [0, \pi],  $$
     for each fixed $n \in \N_0$.  In the  light of  (\ref{cov.matrix.fun.1}), the   covariance function of    $\{ W_n (\varphi), \varphi \in [0, \pi] \}$ is obtained as follows,
          \begin{eqnarray*}
            &    &   \cov (  W_n (\varphi_1),  W_n (\varphi_2))   \\
            & =  &   \int_0^{ 2 \pi}  \int_0^{ 2 \pi} 
                            \cov ( Z (\varphi_1, \theta_1), Z (\varphi_2, \theta_2)  ) \cos (n \theta_1) \cos (n \theta_2) d \theta_1  d \theta_2   \\
       & =  &      \int_0^{ 2 \pi}  \int_0^{ 2 \pi} 
                        C (\varphi_1, \varphi_2, \theta_1-\theta_2) 
                       \cos (n \theta_1) \cos (n \theta_2) d \theta_1  d \theta_2   \\
          & =  &            \int_0^{ 2 \pi}  \int_0^{ 2 \pi}  \sum_{k=0}^\infty   \left\{ b_k (\varphi_1, \varphi_2) \cos (k (\theta_1-\theta_2)) +a_k (\varphi_1, \varphi_2, \theta_1-\theta_2)  \sin (k (\theta_1-\theta_2)) \right\} 
                                 \cos (n \theta_1) \cos (n \theta_2) d \theta_1  d \theta_2   \\
           & =  &     \sum_{k=0}^\infty b_k (\varphi_1, \varphi_2)       \int_0^{ 2 \pi}  \int_0^{ 2 \pi}  
                              ( \cos (k \theta_1) \cos (k \theta_2) +\sin (k \theta_1) \sin (k \theta_2) )
                                 \cos (n \theta_1) \cos (n \theta_2) d \theta_1  d \theta_2   \\
            &    & +  \sum_{k=0}^\infty a_k (\varphi_1, \varphi_2)       \int_0^{ 2 \pi}  \int_0^{ 2 \pi}  
                              ( \sin (k \theta_1) \cos (k \theta_2) -\cos  (k \theta_1) \sin (k \theta_2) )
                                 \cos (n \theta_1) \cos (n \theta_2) d \theta_1  d \theta_2   \\                    
            & =  &     \sum_{k=0}^\infty b_k (\varphi_1, \varphi_2)   \left(
                          \int_0^{ 2 \pi}     \cos ( k \theta_1) \cos (n \theta_1)  d \theta_1   
                           \int_0^{ 2 \pi}     \cos ( k \theta_2) \cos (n \theta_2)  d \theta_2   \right. \\
            &    &    \left.               
                            +  \int_0^{ 2 \pi}     \sin ( k \theta_1) \cos (n \theta_1)  d \theta_1   
                           \int_0^{ 2 \pi}     \sin ( k \theta_2) \cos (n \theta_2)  d \theta_2   \right)      \\
           & =  &
                     \left\{
                     \begin{array}{ll}
                     4 \pi^2 b_0(\varphi_1, \varphi_2),  ~   &   ~  n =0, \\
                          \pi^2 b_n(\varphi_1, \varphi_2),                                                ~   &   ~  n  \in \N,
                      \end{array}        \right.                                                    
                                         \end{eqnarray*}            
           which implies that $b_n(\varphi_1, \varphi_2)$ is a covariance function for each $n \in \N_0$.

\subsection{Proof of Theorem 2}

 The right-hand series of (\ref{series.exp1}) is mean square convergent under the convergent assumption 
 of $\sum\limits_{n=0}^\infty b_n (\varphi_1, \varphi_2)$, since
    \begin{eqnarray*}
    &    & \rE \left|   \sum_{n=n_1}^{n_1+n_2}  ( V_{n1}  ( \varphi) \cos (n \theta) + V_{n2}  ( \varphi) \sin ( n \theta) )      \right|^2 \\
   & =  &   \sum_{n=n_1}^{n_1+n_2} ( E V^2_{n1} (\varphi)   +  E V^2_{n2} (\varphi) )  \\
   & =  &   2 \sum_{n=n_1}^{n_1+n_2} b_n (\varphi, \varphi)   \\
   & \to & 0,   ~~~~~~~~ n_1, n_2 \to \infty.
     \end{eqnarray*}   
  It is easy to verify that     its covariance function  is given by (\ref{cov.matrix.fun.m.0}).

\subsection{Proof of Theorem 3}

  Since $\{ \bZ (\x), \x \in \S^2 \}$ is an $m$-variate axially symmetric and mean square continuous random field, 
 two  scalar random fields $\{ Z_i (\varphi_1, \theta)+Z_j (\varphi_2, \theta), \theta \in [0, 2 \pi] \}$ and $\{ Z_i (\varphi_1, \theta)-Z_j (\varphi_2, \theta), \theta \in [0, 2 \pi] \}$ are also    axially symmetric  
 and mean square continuous for fixed $\varphi_k$ ($k =1, 2$) and $i, j \in \{ 1, \ldots, m \}$, with covariance functions
      \begin{eqnarray*}
     &    &   \cov (  Z_i (\varphi_1, \theta_1 )+Z_j ( \varphi_1, \theta_1 ), Z_i (\varphi_2, \theta_2)+Z_j ( \varphi_2, \theta_2) )  \\
       & = & C_{ii} (\varphi_1, \varphi_2, \theta_1-\theta_2) +C_{ij} ( \varphi_1, \varphi_2, \theta_1-\theta_2) 
                + C_{ji} (\varphi_1, \varphi_2, \theta_1-\theta_2) + C_{jj} (\varphi_1, \varphi_2, \theta_1-\theta_2),  \\
     &  &  \cov (  Z_i (\varphi_1, \theta_1 )-Z_j ( \varphi_1, \theta_1 ), Z_i (\varphi_2, \theta_2)-Z_j ( \varphi_2, \theta_2) )  \\
       & = & C_{ii} (\varphi_1, \varphi_2, \theta_1-\theta_2) -C_{ij} ( \varphi_1, \varphi_2, \theta_1-\theta_2)  
                 -C_{ji} (\varphi_1, \varphi_2, \theta_1-\theta_2) + C_{jj} (\varphi_1, \varphi_2, \theta_1-\theta_2), 
      \end{eqnarray*}
         \hfill   $\varphi_1, \varphi_2,  \theta_1, \theta_2 \in [0, 2 \pi],$
         
         \noindent
         respectively.  It follows from  Theorem 1 that 
 \begin{equation}
 \label{thm2.1}
       \begin{array}{lll}
       &  &   C_{ii} (\varphi_1, \varphi_2, \theta) +C_{ij} ( \varphi_1, \varphi_2, \theta) 
                + C_{ji} (\varphi_1, \varphi_2, \theta) + C_{jj} (\varphi_1, \varphi_2, \theta) \\
       & =  & \sum\limits_{n=0}^\infty  \{ b_{n, ij+} (\varphi_1, \varphi_2) \cos (n \theta) 
                         +a_{n, ij+} (\varphi_1, \varphi_2) \sin (n \theta)  \},  
     \end{array}
 \end{equation}
 and
     \begin{equation}
 \label{thm2.2}
       \begin{array}{lll}
       &  &   C_{ii} (\varphi_1, \varphi_2, \theta) -C_{ij} ( \varphi_1, \varphi_2, \theta) 
                - C_{ji} (\varphi_1, \varphi_2, \theta) + C_{jj} (\varphi_1, \varphi_2, \theta) \\
        & =  & \sum\limits_{n=0}^\infty  \{ b_{n, ij-} (\varphi_1, \varphi_2) \cos (n \theta)
             +a_{n, ij-} (\varphi_1, \varphi_2) \sin (n \theta) \},  
     \end{array}
 \end{equation}
       where $b_{n, ij+} (\varphi_1, \varphi_2)$ and  $b_{n, ij-} (\varphi_1, \varphi_2)$ are covariance functions
       for each $n \in \N_0$, and 
       $\sum\limits_{n=0}^\infty b_{n, ij+} (\varphi_1, \varphi_2)$ and
        $\sum\limits_{n=0}^\infty b_{n, ij-} (\varphi_1, \varphi_2)$ converge.
        Taking the difference between (\ref{thm2.1}) and (\ref{thm2.2}) and then dividing by 4 yields
        $$   \frac{ C_{ij} ( \varphi_1, \varphi_2, \theta) 
                 +  C_{ji} ( \varphi_1, \varphi_2, \theta)      }{2} 
           =  \sum\limits_{n=0}^\infty \{ 
                 b_{n, ij} (\varphi_1, \varphi_2) \cos (n \theta) +a_{n, ij} (\varphi_1, \varphi_2) \sin (n \theta) \},  $$
  where 
  $$ b_{n, ij} = \frac{1}{4} (b_{n, ij+}-b_{n, ij-}), 
      ~~  a_{n, ij} = \frac{1}{4} (a_{n, ij+}-a_{n, ij-}) n \in \N_0. $$
  This establishes (\ref{cov.matrix.fun.m}). 
  
   In order to verify that $\mathbf{B}_n (\varphi_1, \varphi_2)$ is a covariance matrix function
  for each $n \in \N_0$,   consider an $m$-variate stochastic process
     $$ \mathbf{W}_n (\varphi) = \int_0^{ 2 \pi} \frac{ \bZ (\varphi, \theta) +( \tilde{\bZ} (\varphi, \theta) )'}{\sqrt{2}} \cos (n \theta)  d \theta,     ~~~~~~ \varphi \in [0, \pi], $$
     where $\{ \tilde{\bZ} (\x),  \x \in \S^2 \}$ is an independent copy of  $\{ \bZ(\x),  \x \in \S^2 \}$.
   Its covariance matrix function is positively  proportional  to $\mathbf{B}_n (\varphi_1, \varphi_2)$ as follows, 
     \begin{eqnarray*}
     &     &   \cov ( \mathbf{W}_n (\varphi_1), \mathbf{W}_n (\varphi_2))   \\
     & = &  \int_0^{ 2 \pi} \int_0^{ 2 \pi} \frac{ \cov ( \bZ (\varphi_1, \theta_1), \bZ (\varphi_2, \theta_2)) 
                            + \left(  \cov ( \tilde{\bZ}  (\varphi_1, \theta_1),  \tilde{\bZ } (\varphi_2, \theta_2))  \right)' }{2} \cos (n \theta_1) \cos (n \theta_2) d \theta_1  d \theta_2  \\ 
        & =  &    \int_0^{ 2 \pi}  \int_0^{ 2 \pi} 
                       \frac{ \bC (\varphi_1, \varphi_2, \theta_1-\theta_2)
                            +(\bC (\varphi_1, \varphi_2, \theta_1-\theta_2) )' }{2} 
                       \cos (n \theta_1) \cos (n \theta_2) d \theta_1  d \theta_2   \\
               & =  &     \sum_{k=0}^\infty \mathbf{B}_k (\varphi_1, \varphi_2)       \int_0^{ 2 \pi}  \int_0^{ 2 \pi}  
                              ( \cos (k \theta_1) \cos (k \theta_2) +\sin (k \theta_1) \sin (k \theta_2) )
                                 \cos (n \theta_1) \cos (n \theta_2) d \theta_1  d \theta_2   \\
              &    &  +   \sum_{k=0}^\infty \mathbf{A}_k (\varphi_1, \varphi_2)       \int_0^{ 2 \pi}  \int_0^{ 2 \pi}  
                              ( \sin (k \theta_1) \cos (k \theta_2) - \cos (k \theta_1) \sin (k \theta_2) )
                                 \cos (n \theta_1) \cos (n \theta_2) d \theta_1  d \theta_2   \\                 
                  & =  &
                     \left\{
                     \begin{array}{ll}
                     4 \pi^2 \mathbf{B}_0(\varphi_1, \varphi_2),  ~   &   ~  n =0, \\
                          \pi^2 \mathbf{B}_n(\varphi_1, \varphi_2),                                                ~   &   ~  n  \in \N.
                      \end{array}        \right.                                           
     \end{eqnarray*}

  \subsection{Proof of Theorem 5}

    For a fixed $t \in \T$, consider two purely spatial random fields $\left\{  \bZ (\x; 0) + \bZ (\x; t),  \x \in \S^2 \right\}$ and
$\left\{   \bZ (\x; 0) -  \bZ (\x; t),  \x \in \S^2 \right\}$.   Their covariance matrix functions are,  in terms of $\bC( \varphi_1, \varphi_2, \theta; t)$,
     \begin{eqnarray*}
     &   &   \cov \left(     \bZ (\varphi_1, \theta_1; 0) + \bZ (\varphi_1, \theta_1; t), 
                ~ \bZ (\varphi_2, \theta_2; 0) + \bZ (\varphi_2, \theta_2; t)       \right)     \\
     & = &  2  \bC (\varphi_1, \varphi_2, \theta_1- \theta_2; 0) +  \bC ( \varphi_1, \varphi_2, \theta_1- \theta_2; t) 
                +  \bC (\varphi_1, \varphi_2, \theta_1- \theta_2; -t),
     \end{eqnarray*}
     and
       \begin{eqnarray*}
     &   &   \cov \left(     \bZ (\varphi_1, \theta_1; 0) - \bZ (\varphi_1, \theta_1; t),
                ~ \bZ (\varphi_2, \theta_2; 0) - \bZ (\varphi_2, \theta_2; t)       \right)     \\
     & = &  2  \bC (\varphi_1, \varphi_2, \theta_1- \theta_2; 0) - \bC ( \varphi_1, \varphi_2, \theta_1- \theta_2; t) 
                     -  \bC (\varphi_1, \varphi_2, \theta_1- \theta_2; -t),    
                     ~~~ (\varphi_k, \theta_k) \in \S^2,  ~ k=1, 2.
     \end{eqnarray*}
 It follows from Theorem 3  that 
            \begin{equation}
            \label{eq1.thm5}
            \begin{array}{c}
               2  \bC (\varphi_1, \varphi_2, \theta; 0) +  \bC ( \varphi_1, \varphi_2, \theta; t) 
                +  \bC (\varphi_1, \varphi_2, \theta; -t) \\
                + \{ 2  \bC (\varphi_1, \varphi_2, \theta; 0) +  \bC ( \varphi_1, \varphi_2, \theta; t) 
                +  \bC (\varphi_1, \varphi_2, \theta; -t) \}' \\
                 =    \sum\limits_{n=0}^\infty  \left(
                      \mathbf{B}_{n+} (\varphi_1, \varphi_2; t)   \cos (n \theta)
                      + \mathbf{A}_{n+} (\varphi_1, \varphi_2; t)   \sin (n \theta) \right),  
             \end{array}          
           \end{equation}
           and
             \begin{equation}
            \label{eq2.thm5}
              \begin{array}{c}
               2  \bC (\varphi_1, \varphi_2, \theta; 0) -  \bC ( \varphi_1, \varphi_2, \theta; t) 
                -  \bC (\varphi_1, \varphi_2, \theta; -t) \\
                + \{ 2  \bC (\varphi_1, \varphi_2, \theta; 0) - \bC ( \varphi_1, \varphi_2, \theta; t) 
                -  \bC (\varphi_1, \varphi_2, \theta; -t) \}' \\
                 =    \sum\limits_{n=0}^\infty  \left(
                      \mathbf{B}_{n-} (\varphi_1, \varphi_2; t)   \cos (n \theta)
                      + \mathbf{A}_{n-} (\varphi_1, \varphi_2; t)   \sin (n \theta) \right),  
             \end{array}          
     \end{equation}
       \hfill   $ (\varphi_k, \theta_k) \in \S^2,  ~ k=1, 2,  t \in \T.$
       
       \noindent
  Taking the difference between (\ref{eq1.thm5}) and (\ref{eq2.thm5}) results in (\ref{cov.matrix.fun.m.t}), with 
        $$    \mathbf{B}_n (\varphi_1, \varphi_2; t) = \frac{1}{4}  \mathbf{B}_{n+} (\varphi_1, \varphi_2; t) - \frac{1}{4}  \mathbf{B}_{n-} (\varphi_1, \varphi_2; t),  ~~~~~~ $$
      and  $$    \mathbf{A}_n (\varphi_1, \varphi_2; t) = \frac{1}{4}  \mathbf{A}_{n+} (\varphi_1, \varphi_2; t) - \frac{1}{4}  \mathbf{A}_{n-} (\varphi_1, \varphi_2, \theta; t),  ~~~~~~~~   n \in \N_0. $$
        Obviously,   $\mathbf{B}_n(\varphi_1, \varphi_2; t)$ is symmetric,
        and   $\sum\limits_{n=0}^\infty   \mathbf{B}_n(\varphi_1, \varphi_2; t) $ converges. 
        
        In order to verify that $ \mathbf{B}_n (\varphi_1, \varphi_2; t) $ is a  spatio-temporal covariance matrix function on
        $[0, \pi] \times \T$ for each fixed $n \in \N_0$,
        consider an $m$-variate  random field
              $$ \mathbf{W}_n (\varphi; t) = \int_0^{2 \pi}  \frac{ \bZ (\varphi, \theta; t) +\bZ (\varphi, \theta; -t)
                 + ( \tilde{ \bZ } (\varphi, \theta; t) +\tilde{\bZ} (\varphi, \theta; -t) )'}{2} 
                          \cos (n \theta) d \theta,  ~~~ \varphi \in [0, \pi], ~ t \in \T, $$
        where $\{ \tilde{\bZ} (\varphi, \theta; t), \varphi \in [0, \pi], \theta \in [0, 2 \pi],    t \in \T \}$ is an independent copy of 
         $\{ \bZ (\varphi, \theta; t), \varphi \in [0, \pi], \theta \in [0, 2 \pi],    t \in \T \}$.    
            The covariance matrix function of $\{   \mathbf{W}_n (\varphi; t), \varphi \in [0, \pi], t \in \T \}$ is a positive scalar product
            of     $ \mathbf{B}_n (\varphi_1, \varphi_2; t) $, since
             \begin{eqnarray*}
             &   &  \cov (  \mathbf{W}_n (\varphi_1; t_1),  \mathbf{W}_n (\varphi_2; t_2))  \\
             & =  & \frac{1}{2} \int_0^{ 2 \pi}   \int_0^{ 2 \pi} \{ \cov \left(  
                         \bZ (\varphi_1, \theta_1; t_1)+\bZ (\varphi_1, \theta_1; -t_1),
                         ~ \bZ (\varphi_2, \theta_2; t_2)+\bZ (\varphi_2, \theta_2; -t_2) \right)   \\
             &   &    +       
                         \cov (  
                         \tilde{\bZ}'  (\varphi_1, \theta_1; t_1)+ \tilde{\bZ}'   (\varphi_1, \theta_1; -t_1),
                         ~  \tilde{\bZ}'   (\varphi_2, \theta_2; t_2)+ \tilde{\bZ}'   (\varphi_2, \theta_2; -t_2) )  \}     
                         \cos (n \theta_1) \cos ( n \theta_2)  d \theta_1 d \theta_2  \\
             & = &       \frac{1}{2} \int_0^{ 2 \pi}   \int_0^{ 2 \pi} 
                           \{ \bC (\varphi_1, \varphi_2, \theta_1-\theta_2; t_1-t_2)
                                    +     \bC (\varphi_1, \varphi_2, \theta_1-\theta_2; t_2-t_1)  \\
             &    &                          
                                    + (\bC (\varphi_1, \varphi_2, \theta_1-\theta_2; t_1-t_2)
                                    +     \bC (\varphi_1, \varphi_2, \theta_1-\theta_2; t_2-t_1))'    \}
                               \cos (n \theta_1) \cos ( n \theta_2)  d \theta_1 d \theta_2  \\
             & = &             \sum\limits_{n=0}^\infty \int_0^{ 2 \pi}   \int_0^{ 2 \pi}    \left\{   \mathbf{B}_n (\varphi_1, \varphi_2; t) \cos (n  (\theta_1-\theta_2))
                + \mathbf{A}_n (\varphi_1, \varphi_2; t) \sin (n  (\theta_1-\theta_2)) \right\}  \cos (n \theta_1) \cos ( n \theta_2)  d \theta_1 d \theta_2  \\
            &  =  &
                 \left\{
                     \begin{array}{ll}
                     4 \pi^2 \mathbf{B}_0(\varphi_1, \varphi_2; t_1-t_2),  ~   &   ~  n =0, \\
                          \pi^2 \mathbf{B}_n(\varphi_1, \varphi_2; t_1-t_2),                                                ~   &   ~  n  \in \N.
                      \end{array}        \right.
             \end{eqnarray*}

~

\section{Proof of Theorem 6}

 It suffices to prove the last part, while  the proof of the first part is analogous to that of Theorem 2.
 If  all  $\mathbf{V}_{n1} (\varphi; t)$ and $\mathbf{V}_{n2} (\varphi; t)$ ($n \in \N_0$) are normally distributed, then
 $\bZ (\varphi, \theta; t)$ defined by (\ref{series.expm}) is normally distributed, too. On the other hand, if
 $\bZ (\varphi, \theta; t)$ defined by (\ref{series.expm}) is normally distributed, then, for each $n \in \N_0$,
  $\bZ (\varphi, \theta; t) \cos (n \theta)$  and $\bZ (\varphi, \theta; t) \sin (n \theta)$ are normally distributed as well. So are
      \begin{eqnarray*}
       &     &  \int_0^{ 2 \pi}  \bZ (\varphi, \theta; t) \cos (n \theta) d \theta  \\
       &  =  &  \int_0^{ 2 \pi}   \sum_{k=0}^\infty  ( \mathbf{V}_{k1}  ( \varphi; t) \cos (k \theta) + \mathbf{V}_{k2}  ( \varphi; t) \sin ( k \theta) )   \cos (n \theta) d \theta  \\
       &  =  &   \sum_{k=0}^\infty   \mathbf{V}_{k1}  ( \varphi; t)  \int_0^{ 2 \pi}  \cos (k \theta)   \cos (n \theta) d \theta  
                + \sum_{k=0}^\infty   \mathbf{V}_{k2}  ( \varphi; t)  \int_0^{ 2 \pi}  \sin (k \theta)   \cos  (n \theta) d \theta  \\
       &   =   &   \left\{
                      \begin{array}{ll}
                     2 \pi  \mathbf{V}_{01} ( \varphi; t),    ~   &     ~  n =0,  \\
                      \pi  \mathbf{V}_{n1} ( \varphi; t),    ~   &    ~   n \in \N,
                      \end{array}   \right.
      \end{eqnarray*}
      and
      \begin{eqnarray*}
       &     &  \int_0^{ 2 \pi}  \bZ (\varphi, \theta; t) \sin (n \theta) d \theta  \\
        &  =  &   \sum_{k=0}^\infty   \mathbf{V}_{k1}  ( \varphi; t)  \int_0^{ 2 \pi}  \cos (k \theta)   \sin (n \theta) d \theta  
                + \sum_{k=0}^\infty   \mathbf{V}_{k2}  ( \varphi; t)  \int_0^{ 2 \pi}  \sin (k \theta)   \sin  (n \theta) d \theta  \\
       &   =   &   \left\{
                      \begin{array}{ll}
                     0,    ~   &     ~  n =0,  \\
                      \pi  \mathbf{V}_{n2} ( \varphi; t),    ~   &    ~   n \in \N.
                      \end{array}   \right.
      \end{eqnarray*}
      
      ~

      \subsection{Proof of Theorem 7}

  According to Theorem 8 of \cite{Ma2012}, it suffices to check  inequality (\ref{positive.definite}) for  (\ref{cov.m.t.suf}).
  For every $l \in \N$, any $\w_i \in \R^m$, $\varphi_i \in [0, \pi]$,  $\theta_i \in [0, 2 \pi]$, and $t_i \in \T$ ($k=1, \ldots, l$), we have
     \begin{eqnarray*}
     &    &   \sum_{i=1}^l \sum_{j=1}^l \w'_i  \bC(\varphi_i, \varphi_j, \theta_i-\theta_j; t_i -t_j)  \w_j \\
     & =  & \sum_{n=0}^\infty   \sum_{i=1}^l \sum_{j=1}^l \w'_i 
                (\mathbf{B}_n (\varphi_i, \varphi_j; t_i-t_j) \cos (n (\theta_i-\theta_j))
                 +\mathbf{A}_n  (\varphi_i, \varphi_j; t_i-t_j) \cos (n (\theta_i-\theta_j)) )  \w_j \\
     & = &       \sum_{n=0}^\infty   \sum_{i=1}^l \sum_{j=1}^l  \left\{  \w'_i \cos (n \theta_i)  \mathbf{B}_n (\varphi_i, \varphi_j; t_i-t_j)    \w_j  \cos (n \theta_j) 
                  +
                \w'_i  \sin  (n \theta_i) \mathbf{B}_n (\varphi_i, \varphi_j; t_i-t_j)  \w_j  \sin (n \theta_j)    \right.
                  \\
    &     &        \left.     +
                  \w'_i  \sin  (n \theta_i)  \mathbf{A}_n (\varphi_i, \varphi_j; t_i-t_j)  \w_j  \cos (n \theta_j) 
                   - \w'_i  \cos  (n \theta_i) \mathbf{A}_n (\varphi_i, \varphi_j; t_i-t_j)  \w_j  \sin (n \theta_j) 
                  \right\} \\  
     & =  &       \sum_{n=0}^\infty   \sum_{i=1}^l \sum_{j=1}^l  ( \u'_{ni}  \mathbf{B}_n (\varphi_i, \varphi_j; t_i-t_j)   \u_{nj}  
                +    \v'_{ni}  \mathbf{B}_n (\varphi_i, \varphi_j; t_i-t_j)   \v_{nj}  \\
      &   &             + \v'_{ni} \mathbf{A}_n (\varphi_i, \varphi_j; t_i-t_j)  \u_{nj}   - \u'_{ni}  \mathbf{A}_n (\varphi_i, \varphi_j; t_i-t_j)  \v_{nj} 
                    )  \\
     & \ge & 0,           
     \end{eqnarray*}
     where the last inequality follows from the assumption (\ref{suff.con}), $\u_{ni} = \w_i \cos (n \theta_i),
     \v_{ni} = \w_i \sin(n \theta_j). $

~

~


\begin{thebibliography}{10}


 \bibitem{Andrews1999}
  Andrews, G. E.,  Askey, R.,  Roy, R.:
   Special Functions. Cambridge University Press, Cambridge (1999)
   
   
    \bibitem{Askey1976}
 Askey, R., Bingham, N. H.:
 Gaussian processes on compact symmetric spaces.
  Z. Wahrscheinlichkeitstheorie verw. Gebiete {\bf  37},  127-143 (1976)
  
  
  
 
 \bibitem{Bingham1973}
   Bingham, N. H.:
 Positive definite functions on spheres.
   Proc. Cambridge Phil. Soc.  {\bf 73},   145-156 (1973)



\bibitem{Castruccio2013}
Castruccio, S. and Stein, M. L.: Global space-time models for climate ensembles. Ann. Appl. Stat. {\bf 7}, 
  1593-1611 (2013)

  
  \bibitem{Cheng2016}
 Cheng, D.,    Xiao, Y.:
  Excursion  probability of Gaussian random fields on sphere.   Bernoulli {\bf 22},  1113-1130 (2016)
 
 
  \bibitem{Cohen2012}
  Cohen, S.,  Lifshits, M. A.:
  Stationary Gaussian random fields on hyperbolic spaces and on Euclidean spheres.
   ESAIM  {\bf 16},  165-221 (2012)




\bibitem{Dovidio2014}
D'Ovidio, M.:
Coordinates changed random fields on the sphere.
 J. Stat. Phys. {\bf 154}, 1153-1176 (2014)


\bibitem{DuMaLi2013}
Du,  J.,  Ma,  C.,   Li, Y.:
 Isotropic variogram matrix functions on spheres.
 Math.  Geosci. {\bf  45},   341-357 (2013)
 
 
 

\bibitem{Gangolli1967}
 Gangolli, R.:  Positive definite kernels on homogeneous spaces and certain stochastic processes related 
to L{\' e}vy's Brownian motion of several parameters.  Ann Inst H Poincar{\' e} B  {\bf 3},  121-226 (1967)



\bibitem{Gaspari1999}
 Gaspari, G.,    Cohn, S. E.:
 Construction of correlations in two and three dimensions.
  Q. J. R. Meteorol. Soc.  {\bf 125},    723-757 (1999) 

 \bibitem{Gaspari2006}
  Gaspari,  G.,   Cohn,  S. E.,   Guo, J.,   Pawson, S.: 
 Construction and application of covariance functions with variable length-fields.
  Q. J. R. Meteorol. Soc. {\bf 132}, 815-1838 (2006) 
 
 
\bibitem{Gradshteyn2007}
Gradshteyn, I. S.,  Ryzhik,  I. M.:  
Tables of  Integrals, Series, and Products, 7th edtion.
Academic Press, Amsterdam (2007)  

   
\bibitem{Hannan1970}
 Hannan, E. J.:
   Multiple Time  Series.
 Wiley, New York (1970) 
 
 
 \bibitem{Hitczenko2012}
 Hitczenko, M. and Stein, M. L.: Some theory for anisotropic processes on the sphere.
 Statist. Meth. {\bf 9}, 211-227 (2012)
 
 
 \bibitem{Huang2012}
 Huang, C., Zhang, H., and  Robeson, S.: 
 A simplified representation of the covariance structure of axially symmetric processes on the sphere.
 Statist. Probab. Lett. {\bf 82}, 1346-1351 (2012)

 
 \bibitem{Jones1963}
 Jones, R. H.: 
Stochastic processes on a sphere.
 Ann. Math. Statist.  {\bf 34},  213-218 (1963)
 
 \bibitem{Jun2011}
 Jun, M.: Non-stationary cross-covariance models for multivariate processes on a globe.
 Scand. J. Statist.,  {\bf 38},  726-747 ( 2011)
 
 
 
 \bibitem{Jun2007}
 Jun, M. and Stein, M. L.: 
 An approach to producing space-time covariance functions on spheres.
 Technometrics  {\bf 49}, 468-479 (2007) 
 
 \bibitem{Jun2008}
 Jun, M. and Stein, M. L.: 
 Nonstationary covariance models for global data.
  Ann. Appl. Statist.
 {\bf 2},  1271-1289  (2008)
 

\bibitem{Lamberg2001}
Lamberg, L.,    Muinonen,  K.,   Yl$\ddot{\mbox{o}}$nen, J.,     Lumme,  K.:
 Spectral estimation of Gaussian random circles and spheres. 
   J. Comput. Appl. Math.  {\bf 136},   109-121 (2001)
   
   


\bibitem{Leonenko2012}
 Leonenko, N.,    Sakhno, L.:    On spectral representation of tensor random fields on the sphere.
 Stoch. Anal. Appl.  {\bf 31},   167-182 (2012)


\bibitem{Leonenko2013}
Leonenko,  N.,   Shieh,  N.:    R{\' e}nyi function for multifractal random fields.
    Fractals,   {\bf  21},  1350009, 13 pp (2013)
  

 \bibitem{Ma2011}
  Ma, C.:
 Vector random fields with second-order moments or second-order increments.
    Stoch. Anal. Appl.  {\bf 29},   197-215 (2011)
   
 

\bibitem{Ma2012}
  Ma, C.: 
  Stationary and isotropic vector random fields on spheres.
   Math. Geosci.  {\bf  44},   765-778 (2012) 
  

 
 \bibitem{Ma2015}
 Ma,  C.:   Isotropic covariance matrix functions on all  spheres.
 Math. Geosci.  {\bf  47},  699-717 (2015)
 
 \bibitem{Ma2016a}
 Ma, C.:  Stochastic representations of isotropic vector random fields on spheres.
   Stoch. Anal. Appl. {\bf 34}, 389-403 (2016)
   
    \bibitem{Ma2016b}
 Ma, C.:  Isotropic covariance matrix polynomials on spheres.
   Stoch. Anal. Appl.,   To appear. 
   
   
   \bibitem{Ma2016c}
   Ma, C.:  Time varying isotropic vector random fields on spheres.
   J. Theor. Probab., To appear (DOI 10.1007/s10959-016-0689-1).


\bibitem{Malyarenko2013}
  Malyarenko,  A.:
   Invariant Random Fields on Spaces with a Group Action.
 Springer, New York (2013)
    
 \bibitem{Malyarenko1992}
  Malyarenko,  A.,   Olenko, A.:  Multidimensional covariant random fields on commutative locally compact groups.
   Ukrainian Math. J.  {\bf  44},  1384-1389 (1992)

 
  

\bibitem{McLeod1986}
  McLeod,  M. G.:
 Stochastic processes on a sphere.
  Phy. Earth Plan. Interior  {\bf 43},  283-299 (1986)
 
 
\bibitem{Mokljacuk1979}
Mokljacuk, M. P., Jadrenko, M. I.:
Linear statistical problems for stationary isotropic random fields on a sphere, I.
Theor. Prob. Math. Statist. {\bf 18}, 115-124 (1979)


 
 
 
 \bibitem{Roy1973}
  Roy, R.:  
  Spetral analysis for random process on the circle.
  J. Appl. Prob. {\bf 9}, 745-757 (1972)
 
  
 \bibitem{Roy1976}
   Roy,  R.:
 Spectral analysis for a random process on the sphere.
 Ann. Inst. Statist. Math.  {\bf 28},   91-97 (1976) 



 \bibitem{Schoenberg1942}
  Schoenberg, I.:   Positive definite functions on spheres.
  Duke Math. J.  {\bf 9},   96-108 (1942)
  
  
  \bibitem{Stein2007}
  Stein, M.: Spatial variation of total column ozone on a global scale.
  Ann. Appl. Statist. {\bf 1}, 191-210 (2007)

 \bibitem{Szego1959}
  Szeg$\ddot{\mbox{o}}$, G. 
  Orthogonal Polynomials,  4th edition. Amer. Math. Soc. Colloq. Publ., vol 23. Amer. Math. Soc.,
Providence (1975)

 \bibitem{Watson1944}
   Watson,  G. N.:  A  Treatise on the Theory of Bessel Functions, 2nd edn. Cambridge University Press,
London (1944)

\bibitem{Yadrenko1983}
  Yadrenko,  A. M. 
  Spectral  Theory of  Random Fields.
Optimization Software, New York (1983)

\bibitem{Yaglom1961}
 Yaglom,  A. M.: 
Second-order homogeneous random fields. Proc. 4th Berkeley Symp. Math. Stat. Prob.  {\bf 2},   593-622 (1961)

\bibitem{Yaglom1987}
  Yaglom,  A. M.:
  Correlation  Theory of  Stationary and  Related  Random  Functions.
 vol. I.  Springer, New York (1987)
 
 
 


 
 
\end{thebibliography}
\end{document}